\documentclass[12pt, reqno]{amsart}
\usepackage{amsmath, amsthm, amscd, amsfonts, amssymb, graphicx, color}
\usepackage[bookmarksnumbered, colorlinks, plainpages]{hyperref}

\hypersetup{colorlinks=true,linkcolor=red, anchorcolor=green, citecolor=cyan, urlcolor=red, filecolor=magenta, pdftoolbar=true}
\newtheorem{theorem}{Theorem}[section]

\theoremstyle{definition}

\theoremstyle{remark}
\newtheorem{remark}[theorem]{Remark}
\numberwithin{equation}{section}
\input{tcilatex}

\begin{document}
\title[Existence of solutions for a class of quasilinear systems]{\textbf{On
the} \textbf{radial solutions of a stationary Schr\"{o}dinger system with a
nonlinear random operator }}
\author[D.-P. Covei]{\textsf{Dragos-Patru Covei}$^{1}$\textsf{\ }}
\address{ $^{1}${\small \textit{Department of Applied Mathematics}}\\
{\small \ \textit{The Bucharest University of Economic Studies }}\\
{\small \textit{Piata Romana, 1st district, postal code: 010374, postal
office: 22, Romania}}}
\email{\textcolor[rgb]{0.00,0.00,0.84}{coveid@yahoo.com}}
\subjclass[2010]{Primary: 35J55, 35J60 Secondary: 35J65.}
\keywords{Entire solution; Fixed Point; Large solution; Bounded Solution;
Nonlinear system.}
\date{Received: xxxxxx; Revised: yyyyyy; Accepted: zzzzzz. \\
\indent $^{*}$Corresponding author}

\begin{abstract}
Our main goal is to establish sufficient conditions for the existence of
positive entire radially symmetric solutions for a system of the type 
\begin{equation*}
\left\{ 
\begin{array}{l}
\func{div}(\phi _{1}(|\nabla u|)\nabla u)+\sigma _{1}\left( \left\vert
x\right\vert \right) \phi _{1}(|\nabla u|)\left\vert \nabla u\right\vert
=p_{1}(\left\vert x\right\vert )f_{1}(u,v),\text{ }x\in \mathbb{R}^{N}\text{ 
}, \\ 
\func{div}(\phi _{2}(|\nabla v|)\nabla v)+\sigma _{2}\left( \left\vert
x\right\vert \right) \phi _{2}(|\nabla v|)\left\vert \nabla v\right\vert
=p_{2}(\left\vert x\right\vert )f_{2}(u,v),\text{ }x\in \mathbb{R}^{N}\text{ 
},%
\end{array}%
\right.
\end{equation*}%
where $\phi _{1},$ $\phi _{2},\ a_{1},$ $a_{2},$ $p_{1}$, $p_{2}$, $f_{1}$
and $f_{2}$ are continuous functions satisfying certain properties. Our
results are obtained by an application of the Arzela--Ascoli theorem.
\end{abstract}

\maketitle

\setcounter{page}{1}


\let\thefootnote\relax\footnote{%
Copyright 2016 by the JOURNAL.}

\section{INTRODUCTION}

The existence of positive entire solutions for the coupled nonlinear systems 
\begin{equation}
\left\{ 
\begin{array}{l}
\Delta _{\phi _{1}}u+\sigma _{1}\left( \left\vert x\right\vert \right) \phi
_{1}(|\nabla u|)\left\vert \nabla u\right\vert =p_{1}(\left\vert
x\right\vert )f_{1}(u,v)\text{ for }x\in \mathbb{R}^{N}, \\ 
\Delta _{\phi _{2}}v+\sigma _{2}\left( \left\vert x\right\vert \right) \phi
_{2}(|\nabla v|)\left\vert \nabla v\right\vert =p_{2}(\left\vert
x\right\vert )f_{2}(u,v)\text{ for }x\in \mathbb{R}^{N}\text{ },%
\end{array}%
\right.  \label{11}
\end{equation}%
where $\mathbb{R}^{N}$ ($N\geq 3$) denote the Euclidean $N$-space, $%
\left\vert \circ \right\vert $ will denote any $N$-dimensional norm, $\Delta
_{\phi _{i}}w$ $(i=1,2)$ stands for the $\phi _{i}$-Laplacian operator
defined as $\Delta _{\phi _{i}}w:=\func{div}(\phi _{i}(|\nabla w|)\nabla w)$
and the functions $\phi _{i}$ satisfy, throughout this paper:

O1)\quad $\phi _{i}$ $\in C^{1}\left( \left( 0,\infty \right) ,\left(
0,\infty \right) \right) $,

O2)\quad $t\phi _{i}(t)$ is strictly increasing in $\left( 0,\infty \right) $%
,

O3)\quad there exist $l_{i}$, $m_{i}>1$ such that 
\begin{equation*}
\text{if }\Phi _{i}\left( t\right) :=\int_{0}^{t}s\phi _{i}\left( s\right) ds%
\text{ then }l_{i}\leq \frac{\Phi _{i}^{\prime }\left( t\right) \cdot t}{%
\Phi _{i}\left( t\right) }\leq m_{i}\text{ for any }t>0,
\end{equation*}

O4)\quad there exist $a_{0}^{i}$, $a_{1}^{i}>0$ such that 
\begin{equation*}
a_{0}^{i}\leq \frac{\Phi _{i}^{\prime \prime }\left( t\right) \cdot t}{\Phi
_{i}^{\prime }\left( t\right) }\leq a_{1}^{i}\text{ for any }t>0\text{,}
\end{equation*}%
have been intensively studied in the last few decades in view of the
understanding of some basic phenomena arising in physics (for more details,
see \cite{DP,CD2}, Kawano-Kusano \cite{KK}, Franchi-Lanconelli-Serrin \cite%
{FLS}, Fukagai-Narukawa \cite{FK}, Grosse-Martin \cite{GR} and Smooke \cite%
{S}). Below there are some examples of functions $\varphi _{1}$ and $\phi
_{2}$ that fulfil (O1)-(O4) which, cf. \cite{FK}, arise in mathematical
models in nonlinear physical science:

\textbf{E1:\quad Nonlinear Elasticity:} 
\begin{equation*}
\Phi _{i}\left( t\right) =\left( 1+t^{2}\right) ^{p}-1,\text{ }\phi
_{i}\left( t\right) =2p\left( 1+t^{2}\right) ^{p-1},
\end{equation*}%
where $t>0$ and $p>\frac{1}{2}$;

\textbf{E2:\quad Plasticity:} 
\begin{equation*}
\Phi _{i}\left( t\right) =t^{p}\left( \ln \left( 1+t\right) \right) ^{q},%
\text{ }\phi _{i}\left( t\right) =\frac{\ln ^{q-1}\left( t+1\right) }{t+1}%
\left[ \left( pt^{p-1}+pt^{p-2}\right) \ln \left( t+1\right) +qt^{p-1}\right]
\allowbreak ,
\end{equation*}%
where $t>0$, $p>1$ and $q>0$;

\textbf{E3:\quad Generalized Newtonian fluids:} 
\begin{equation*}
\Phi _{i}\left( t\right) =\int_{0}^{t}s^{1-p}\left( \sinh ^{-1}s\right)
^{q}ds,\text{ }\phi _{i}\left( t\right) =\allowbreak t^{-p}\func{arcsinh}%
^{q}t,
\end{equation*}%
where $t>0$, $0\leq p\leq 1$ and $q>0$;

\textbf{E4:}\quad \textbf{Plasma Physics:} 
\begin{equation*}
\Phi _{i}\left( t\right) =\frac{t^{p}}{p}+\frac{t^{q}}{q},\phi _{i}\left(
t\right) =t^{p-2}+t^{q-2},
\end{equation*}%
where $t>0$ and $1<p<q$.

\textbf{E5:\quad Non-Newtonian Fluid:}%
\begin{equation*}
\Phi _{i}\left( t\right) =\frac{t^{p}}{p},\phi _{i}\left( t\right) =t^{p-2},
\end{equation*}%
where $t>0$ and $p>1$.

\begin{remark}
The systems of the form (\ref{11}) are known today as coupled nonlinear
systems of Schr\"{o}dinger type. Since, in particularly, one of the most
important classes of (\ref{11}) is the time-independent Schr\"{o}dinger
equation in quantum mechanics%
\begin{equation}
\Delta u=\frac{8\pi ^{2}m}{h^{2}}\left( V\left( r\right) -E\right) u\text{,}
\label{sc}
\end{equation}%
where $m$ is the particle's "reduced mass", $V\left( r\right) $ is its
potential energy, $r$ is a position vector, $h$ is the Planck constant, $E$
is the energy of a photon and the unknown function $u$ is the wave function
(for more details, see \cite{GR}, \cite{L}, \cite{SJ}, \cite{S} and \cite{XZ}%
).
\end{remark}

In the literature, an \textit{entire large solution} means a couple $\left(
u,v\right) \in C^{1}(\left[ 0,\infty \right) )\times C^{1}(\left[ 0,\infty
\right) )$ of positive functions satisfying (\ref{11}) and such that both $%
u\left( x\right) $ and $v\left( x\right) $ tend to infinity as $\left\vert
x\right\vert \rightarrow \infty $; an \textit{entire} \textit{bounded
solution} if the condition $u\left( x\right) <\infty $ and $v\left( x\right)
<\infty $ as $\left\vert x\right\vert \rightarrow \infty $; a \textit{%
semifinite entire large solution} when ($u\left( x\right) <\infty $ and $%
v\left( x\right) $ tend to infinity) or ($u\left( x\right) $ tend to
infinity and $v\left( x\right) <\infty $) as $\left\vert x\right\vert
\rightarrow \infty $.

In the next, we shall reserve $r$ for the polar distance, $r:=\sqrt{%
x_{1}^{2}+...+x_{n}^{2}}$ for $x=\left( x_{1},...,x_{n}\right) \in \mathbb{R}%
^{N}$. Note that if $\phi _{i}$ is considered as a function in $\mathbb{R}%
^{N}$ depending only on $r$, then 
\begin{equation*}
\Delta _{\phi _{i}}w=\left( \phi _{i}(\left\vert w\right\vert ^{\prime
})\left\vert w\right\vert ^{\prime }\right) ^{\prime }+\frac{N-1}{r}\phi
(\left\vert w\right\vert ^{\prime })\left\vert w\right\vert ^{\prime }.
\end{equation*}%
We would like to quote some references where the existence of entire bounded
radial solutions, or the existence of entire large radial solutions for the
systems of the form (\ref{11}) were analyzed. Lair, \cite{L1} has considered
the entire large radial solutions for the elliptic system 
\begin{equation}
\left\{ 
\begin{array}{l}
\Delta u=p_{1}\left( r\right) v^{\alpha }, \\ 
\Delta v=p_{2}\left( r\right) u^{\beta }\text{, }r=\left\vert x\right\vert 
\text{, }x\in \mathbb{R}^{N}\text{ (}N\geq 3\text{),}%
\end{array}%
\right.  \label{lair}
\end{equation}%
where $0<\alpha \leq 1$, $0<\beta \leq 1$, $p_{1}$ and $p_{2}$ are
nonnegative continuous functions on $\mathbb{R}^{N}$. He proved that a
necessary and sufficient condition for this system to have a positive entire
large radial solution, is 
\begin{eqnarray}
\int_{0}^{\infty }tp_{1}\left( t\right) \left(
t^{2-N}\int_{0}^{t}s^{N-3}Q\left( s\right) ds\right) ^{\alpha }dt &=&\infty ,
\label{c1l} \\
\int_{0}^{\infty }tp_{2}\left( t\right) \left(
t^{2-N}\int_{0}^{t}s^{N-3}P\left( s\right) ds\right) ^{\beta }dt &=&\infty ,
\label{c2l}
\end{eqnarray}%
where $P\left( r\right) =\int_{0}^{r}\tau p_{1}\left( \tau \right) d\tau $
and $Q\left( r\right) =\int_{0}^{r}\tau p_{2}\left( \tau \right) d\tau $.

It is well known, see Yang \cite{YH}, that if $p:\left[ 0,\infty \right)
\rightarrow \left[ 0,\infty \right) $ is a spherically symmetric continuous
function and the nonlinearity $f:[0,\infty )\rightarrow \lbrack 0,\infty )$
is a continuous, increasing function with $f\left( 0\right) \geq 0$ and $%
f\left( s\right) >0$ for all $s>0$ which satisfies%
\begin{equation}
\int_{1}^{\infty }\frac{1}{f\left( t\right) }dt=\infty ,  \label{DY}
\end{equation}%
then, the single equation%
\begin{equation}
\Delta u=p\left( r\right) f\left( u\right) \text{ for }x\in \mathbb{R}^{N}%
\text{ (}N\geq 3\text{), }\underset{r\rightarrow \infty }{\lim }u\left(
r\right) =\infty \text{,}  \label{dye}
\end{equation}%
has a nonnegative radial solution if and only if $p$ satisfies 
\begin{equation*}
\underset{r\rightarrow \infty }{\lim }\mathcal{P}_{p}\left( r\right) =\infty 
\text{, }\mathcal{P}_{p}\left( r\right)
:=\int_{0}^{r}s^{1-N}\int_{0}^{s}z^{N-1}p(z)dzds.
\end{equation*}%
A direct computation gives 
\begin{equation*}
\underset{r\rightarrow \infty }{\lim }\mathcal{P}_{p}\left( r\right) =\frac{1%
}{N-2}\int_{0}^{\infty }tp\left( t\right) dt.
\end{equation*}%
However, there is no existence results for the system (\ref{11}) where $%
f_{1} $ and $f_{2}$ satisfy a condition of the form (\ref{DY}). This
observation can be found in the paper of \cite{L1}. Fang-Yi \cite{BF}, in
the particular case $\varphi _{1}\left( t\right) =\varphi _{2}\left(
t\right) =t^{p-1}$ ($p>1$), supplied a sufficient condition 
\begin{equation}
\int_{a}^{\infty }\frac{1}{f_{1}^{1/(p-1)}\left( t,t\right)
+f_{2}^{1/(p-1)}\left( t,t\right) }dt=\infty ,\text{ }t\geq a>0,  \label{zz}
\end{equation}%
for the existence of positive radial large solutions to (\ref{11}). The
condition (\ref{zz}) have been used by many authors and in many contexts,
see Li-Zhang-Zhang \cite{ZL}, Liu-Zhang \cite{HL}, Qin-Yang \cite{QH},
Dkhil-Zeddini \cite{DZ}, and the references therein.

Now we return to (\ref{zz}): if the function $%
(f_{1}^{1/(p-1)}+f_{2}^{1/(p-1)})$ satisfies condition (\ref{zz}), so do
separately each of the functions $f_{1}^{1/(p-1)}$ and $f_{2}^{1/(p-1)}$ but
the converse is not true as one can see from the paper of Bernfeld \cite[%
Example 3.8., pp. 283]{B}. One of our main purposes of this paper is to
establish sufficient conditions for the existence of entire large radial
solutions of the system (\ref{11}) under the new conditions of the form%
\begin{equation}
\int_{a}^{\infty }\frac{1}{f_{1}^{1/(p-1)}\left( t,t\right) }dt=\infty ,%
\text{ }t\geq a>0,  \label{cov1}
\end{equation}%
and%
\begin{equation}
\int_{a}^{\infty }\frac{1}{f_{2}^{1/(p-1)}\left( t,t\right) }dt=\infty ,%
\text{ }t\geq a>0,  \label{cov2}
\end{equation}%
The existence of entire bounded/semifinite entire large positive solutions
is also studied in this paper. Finally, we should like to mention that the
method presented here also yields much more precise information on the
behavior of solutions.

\section{NOTATIONS AND PRELIMINARIES}

We work under the following assumptions:

(P1)\quad $\sigma _{1},\sigma _{2},p_{1},p_{2}:\left[ 0,\infty \right)
\rightarrow \left[ 0,\infty \right) $ are continuous functions\textit{;}

(C1)\quad $f_{1},f_{2}:\left[ 0,\infty \right) \times \left[ 0,\infty
\right) \rightarrow \left[ 0,\infty \right) $ are continuous, increasing, $%
f_{1}\left( 0,0\right) \geq 0$, $f_{2}\left( 0,0\right) \geq 0$ and $%
f_{1}\left( s_{1},s_{2}\right) >0$, $f_{2}\left( s_{1},s_{2}\right) >0$
whenever $s_{1},s_{2}>0$;

(C2)\quad there exist the continuous and increasing functions $h_{1}$, $%
h_{2}:\left[ 0,\infty \right) \times \left[ 0,\infty \right) \rightarrow %
\left[ 0,\infty \right) $ and $\overline{f}_{1}$, $\overline{f}_{2}:\left[
0,\infty \right) \rightarrow \left[ 0,\infty \right) $ such that 
\begin{eqnarray}
f_{1}\left( t_{1},t_{1}\cdot s_{1}\right) &\leq &h_{1}\left(
t_{1},t_{1}\right) \cdot \overline{f}_{1}\left( s_{1}\right) \text{, }%
\forall s_{1}\geq 1\text{ and }\forall \text{ }t_{1}\geq M_{1}a_{1},
\label{c22} \\
f_{2}\left( t_{2},t_{2}\cdot s_{2}\right) &\leq &h_{2}\left(
t_{2},t_{2}\right) \cdot \overline{f}_{2}\left( s_{2}\right) \text{, }%
\forall s_{2}\geq 1\text{ and }\forall \text{ }t_{2}\geq M_{2}a_{2},
\label{c222}
\end{eqnarray}%
where $a_{1},a_{2}\in \left( 0,\infty \right) $, $M_{1}\geq \max \left\{ 1,%
\frac{1}{a_{1}}\right\} $ and $M_{2}\geq \max \left\{ 1,\frac{1}{a_{2}}%
\right\} $.

In order to state our existence theorems, we introduce the following
notations%
\begin{eqnarray*}
Z\left( r\right) &=&\int_{a_{1}+a_{2}}^{r}\frac{1}{\overline{\theta }%
_{1}\left( f_{1}\left( t,t\right) \right) +\overline{\theta }_{2}\left(
f_{2}\left( t,t\right) \right) }dt\text{, }\mathcal{H}_{i}\left( r\right)
=\int_{a_{i}}^{r}\frac{1}{\overline{\theta }_{i}\left( h_{i}\left(
t,M_{i}t\right) \right) }dt\text{, } \\
\xi _{i}\left( t\right) &=&t^{N-1}e^{\int_{0}^{t}\sigma _{i}\left( s\right)
ds}\text{, }P_{i}\left( r\right) =\int_{0}^{r}\Psi _{i}^{-1}\left( \frac{1}{%
\xi _{i}\left( z\right) }\int_{0}^{z}\xi _{i}\left( t\right) p_{i}\left(
t\right) dt\right) dz\text{,} \\
\overline{P}_{i}\left( r\right) &=&\int_{0}^{r}\Psi _{i}^{-1}\left( \frac{1}{%
\xi _{i}\left( t\right) }\int_{0}^{t}\xi _{i}\left( s\right) p_{i}\left(
s\right) \overline{f}_{i}\left( 1+Z^{-1}\left( P_{1}\left( s\right)
+P_{2}\left( s\right) \right) \right) ds\right) dt\text{, } \\
\underline{P}_{1}\left( r\right) &=&\int_{0}^{r}\Psi _{1}^{-1}\left( \frac{1%
}{\xi _{1}\left( t\right) }\int_{0}^{t}\xi _{1}\left( s\right) p_{1}\left(
s\right) f_{1}\left( a_{1},a_{2}+\underline{\theta }_{2}(f_{2}\left(
a_{1},a_{2}\right) )P_{2}\left( s\right) \right) ds\right) dt\text{, } \\
\underline{P}_{2}\left( r\right) &=&\int_{0}^{r}\Psi _{2}^{-1}\left( \frac{1%
}{\xi _{2}\left( t\right) }\int_{0}^{t}\xi _{2}\left( s\right) p_{2}\left(
s\right) f_{2}\left( a_{1}+\underline{\theta }_{1}(f_{1}\left(
a_{1},a_{2}\right) P_{1}\left( s\right) ,a_{2}\right) ds\right) dt\text{, }
\\
\mathcal{H}_{i}\left( \infty \right) &=&\lim_{s\rightarrow \infty }\mathcal{H%
}_{i}\left( s\right) \text{, }P_{i}\left( \infty \right) =\lim_{r\rightarrow
\infty }P_{i}\left( r\right) \text{, \ }i=1,2,\text{ } \\
\overline{P}_{i}\left( \infty \right) &=&\lim_{r\rightarrow \infty }%
\overline{P}_{i}\left( r\right) \text{, }\underline{P}_{i}\left( \infty
\right) =\lim_{r\rightarrow \infty }\underline{P}_{i}\left( r\right) \text{. 
}
\end{eqnarray*}%
Some remarks are now in order on the preliminaries stated in this and the
previous section.

\begin{remark}
A simple example of $f_{1}$ and $f_{2}$ satisfying (C2) is given by \ $%
f_{1}\left( u,v\right) =u^{\beta _{1}}v^{\alpha _{1}}$ and $f_{2}\left(
u,v\right) =u^{\beta _{2}}v^{\alpha _{2}}$ with $\alpha _{1}$, $\beta _{1}$, 
$\alpha _{2}$, $\beta _{2}\in \left[ 0,\infty \right) $ with $\alpha
_{1}^{2}+$ $\beta _{1}^{2}\neq 0$ and $\alpha _{2}^{2}+\beta _{2}^{2}\neq 0$.
\end{remark}

For the proof of the next remark that will be stated below, we refer the
reader to \cite[Lemma 2.1]{FK}.

\begin{remark}
Suppose $\phi _{i}$ ($i=1,2$) satisfy (O1), (O2), (O3) and (O4). Then,%
\begin{equation}
\underline{\theta }_{i}(s_{1})\Psi _{i}^{-1}(s_{2})\leq \Psi
_{i}^{-1}(s_{1}s_{2})\leq \overline{\theta }_{i}(s_{1})\Psi _{i}^{-1}(s_{2})%
\text{ for all }s_{1},s_{2}>0,  \label{ineq}
\end{equation}%
where $\underline{\theta }_{i}(t)=\min \left\{
t^{1/m_{i}},t^{1/l_{i}}\right\} $, $\overline{\theta }_{i}(t)=\max \left\{
t^{1/m_{i}},t^{1/l_{i}}\right\} $.
\end{remark}

The reader is referred to Krasnosel'skii and Rutickii \cite{KR} (see also
Rao and Ren \cite{RAO}) for a through treatment of the assumptions (C2) and (%
\ref{ineq}).

\begin{remark}
If $P_{i}\left( \infty \right) =\infty $ then $\underline{P}_{i}\left(
\infty \right) =\infty $ and $\overline{P}_{i}\left( \infty \right) =\infty $%
. On the other other hand, if $\underline{P}_{i}\left( \infty \right)
=\infty $ or $\overline{P}_{i}\left( \infty \right) =\infty $ then we can
have one of the following 
\begin{equation*}
\begin{array}{cc}
1. & P_{1}\left( \infty \right) <\infty \text{ and }P_{2}\left( \infty
\right) =\infty , \\ 
2. & P_{1}\left( \infty \right) =\infty \text{ and }P_{2}\left( \infty
\right) <\infty , \\ 
3. & P_{1}\left( \infty \right) =\infty \text{ and }P_{2}\left( \infty
\right) =\infty ,%
\end{array}%
\end{equation*}%
(see \cite{L1} for an example in this direction).
\end{remark}

\section{STATEMENTS AND PROOFS OF THE THEOREMS}

Our main objective in this work is to prove the following result:

\begin{theorem}
\label{th1}The system (\ref{11}) has one positive radial solution $\left(
u,v\right) \in C^{1}\left( \left[ 0,\infty \right) \right) \times
C^{1}\left( \left[ 0,\infty \right) \right) $ given \textit{that }$\mathcal{H%
}_{1}\left( \infty \right) =\mathcal{H}_{2}\left( \infty \right) =\infty $
and \textrm{(P1)}, \textrm{(C1)}, \textrm{(C2) }hold true. Moreover, if $%
\underline{P}_{1}\left( \infty \right) =\infty $ and $\underline{P}%
_{2}\left( \infty \right) =\infty $ then 
\begin{equation*}
\lim_{r\rightarrow \infty }u\left( r\right) =\infty \text{ and }%
\lim_{r\rightarrow \infty }v\left( r\right) =\infty .
\end{equation*}
\end{theorem}

\subparagraph{\textbf{Proof of Theorem \protect\ref{th1}:}}

We start by showing that (\ref{11}) has positive radial solutions. On this
purpose we can see that radial solutions of the system 
\begin{equation}
\left\{ 
\begin{array}{l}
\left( \phi _{1}(u^{\prime })u^{\prime }\right) ^{\prime }+\frac{N-1}{r}\phi
_{1}(u^{\prime })u^{\prime }+\sigma _{1}\left( r\right) \phi _{1}(u^{\prime
})u^{\prime }=p_{1}\left( r\right) f_{1}\left( u\left( r\right) ,v\left(
r\right) \right) \text{, }r\geq 0, \\ 
\left( \phi _{2}(v^{\prime })v^{\prime }\right) ^{\prime }+\frac{N-1}{r}\phi
_{2}(v^{\prime })v^{\prime }+\sigma _{2}\left( r\right) \phi _{2}(v^{\prime
})v^{\prime }=p_{2}\left( r\right) f_{2}\left( u\left( r\right) ,v\left(
r\right) \right) \text{, }r\geq 0, \\ 
u^{\prime },v^{\prime }\geq 0\text{ on }\left[ 0,\infty \right)  \\ 
u\left( 0\right) =a_{1}\text{, }v\left( 0\right) =a_{2}%
\end{array}%
\right.   \label{ss1}
\end{equation}%
solve (\ref{11}). By the symmetry of $\left( u,v\right) $ and using the
standard integrating procedure, we rewrite the system (\ref{ss1}) as 
\begin{equation}
\left\{ 
\begin{array}{l}
u\left( r\right) =a_{1}+\int_{0}^{r}\Psi _{1}^{-1}\left( \frac{1}{\xi
_{1}\left( t\right) }\int_{0}^{t}\xi _{1}\left( s\right) p_{1}\left(
s\right) f_{1}\left( u\left( s\right) ,v\left( s\right) \right) ds\right) dt%
\text{, }r\geq 0, \\ 
v\left( r\right) =a_{2}+\int_{0}^{r}\Psi _{2}^{-1}\left( \frac{1}{\xi
_{2}\left( t\right) }\int_{0}^{t}\xi _{2}\left( s\right) p_{2}\left(
s\right) f_{2}\left( u\left( s\right) ,v\left( s\right) \right) ds\right) dt%
\text{, }r\geq 0.%
\end{array}%
\right.   \label{ss}
\end{equation}%
The solution $\left( u,v\right) $ can be constructed by the following
approximate scheme: define $u_{0}=a_{1},v_{0}=a_{2}$ and let $\left\{ \left(
u_{n},v_{n}\right) \right\} _{n\geq 1}$ on $\left[ 0,\infty \right) \times %
\left[ 0,\infty \right) $ given by 
\begin{equation}
\left\{ 
\begin{array}{l}
u_{n}\left( r\right) =a_{1}+\int_{0}^{r}\Psi _{1}^{-1}\left( \frac{1}{\xi
_{1}\left( t\right) }\int_{0}^{t}\xi _{1}\left( s\right) p_{1}\left(
s\right) f_{1}\left( u_{n-1}\left( s\right) ,v_{n-1}\left( s\right) \right)
ds\right) dt\text{, }r\geq 0, \\ 
v_{n}\left( r\right) =a_{2}+\int_{0}^{r}\Psi _{2}^{-1}\left( \frac{1}{\xi
_{2}\left( t\right) }\int_{0}^{t}\xi _{2}\left( s\right) p_{2}\left(
s\right) f_{2}\left( u_{n-1}\left( s\right) ,v_{n-1}\left( s\right) \right)
ds\right) dt\text{, }r\geq 0.%
\end{array}%
\right.   \label{recs}
\end{equation}%
We show that $\left\{ u_{n}\right\} _{n\geq 0}$ and $\left\{ v_{n}\right\}
_{n\geq 0}$ are nondecreasing on $\left[ 0,\infty \right) $. To see this,
express 
\begin{eqnarray*}
u_{1}\left( r\right)  &=&a_{1}+\int_{0}^{r}\Psi _{1}^{-1}\left( \frac{1}{\xi
_{1}\left( t\right) }\int_{0}^{t}\xi _{1}\left( s\right) p_{1}\left(
s\right) f_{1}\left( u_{0}\left( s\right) ,v_{0}\left( s\right) \right)
ds\right) dt \\
&=&a_{1}+\int_{0}^{r}\Psi _{1}^{-1}\left( \frac{1}{\xi _{1}\left( t\right) }%
\int_{0}^{t}\xi _{1}\left( s\right) p_{1}\left( s\right) f_{1}\left(
a_{1},a_{2}\right) ds\right) dt \\
&\leq &a_{1}+\int_{0}^{r}\Psi _{1}^{-1}\left( \frac{1}{\xi _{1}\left(
t\right) }\int_{0}^{t}\xi _{1}\left( s\right) p_{1}\left( s\right)
f_{1}\left( u_{1}\left( s\right) ,v_{1}\left( s\right) \right) ds\right)
dt=u_{2}\left( r\right) .
\end{eqnarray*}%
This proves that $u_{1}\left( r\right) \leq u_{2}\left( r\right) $.
Similarly, $v_{1}\left( r\right) \leq v_{2}\left( r\right) $. By an
induction argument we get 
\begin{equation*}
u_{n}\left( r\right) \leq u_{n+1}\left( r\right) \text{ for any }n\in 
\mathbb{N}\text{ and }r\in \left[ 0,\infty \right) ,
\end{equation*}%
and%
\begin{equation*}
v_{n}\left( r\right) \leq v_{n+1}\left( r\right) \text{ for any }n\in 
\mathbb{N}\text{ and }r\in \left[ 0,\infty \right) .
\end{equation*}%
Let us now prove that the non-decreasing sequences $\left\{ u_{n}\right\}
_{n\geq 0}$ and $\left\{ v_{n}\right\} _{n\geq 0}$ are bounded from above on
bounded sets. By the monotonicity of $\left\{ u_{n}\right\} _{n\geq 0}$ and $%
\left\{ v_{n}\right\} _{n\geq 0}$ one get%
\begin{eqnarray}
\left[ \xi _{1}\left( r\right) \phi _{1}(u_{n}^{\prime }\left( r\right)
)u_{n}^{\prime }\left( r\right) \right] ^{\prime } &=&\xi _{1}\left(
r\right) p_{1}\left( r\right) f_{1}\left( u_{n-1}\left( r\right)
,v_{n-1}\left( r\right) \right)   \notag \\
&\leq &\xi _{1}\left( r\right) p_{1}\left( r\right) f_{1}\left( u_{n}\left(
r\right) ,v_{n}\left( r\right) \right) ,  \label{gen1} \\
\left[ \xi _{2}\left( r\right) \phi _{2}(v_{n}^{\prime }\left( r\right)
)v_{n}^{\prime }\left( r\right) \right] ^{\prime } &\leq &\xi _{2}\left(
r\right) p_{2}\left( r\right) f_{2}\left( u_{n}\left( r\right) ,v_{n}\left(
r\right) \right) .  \label{gen2}
\end{eqnarray}%
Integrating the above inequalities and using (\ref{ineq}), yield that%
\begin{eqnarray*}
u_{n}^{\prime }\left( r\right)  &\leq &\Psi _{1}^{-1}\left( \frac{1}{\xi
_{1}\left( r\right) }\int_{0}^{r}\xi _{1}\left( s\right) p_{1}\left(
s\right) f_{1}\left( u_{n}\left( s\right) ,v_{n}\left( s\right) \right)
ds\right)  \\
&\leq &\Psi _{1}^{-1}\left( \frac{f_{1}\left( u_{n}\left( r\right)
,v_{n}\left( r\right) \right) }{\xi _{1}\left( r\right) }\int_{0}^{r}\xi
_{1}\left( s\right) p_{1}\left( s\right) ds\right)  \\
&\leq &\overline{\theta }_{1}\left( f_{1}\left( u_{n}\left( r\right)
,v_{n}\left( r\right) \right) \right) \Psi _{1}^{-1}\left( \frac{1}{\xi
_{1}\left( r\right) }\int_{0}^{r}\xi _{1}\left( s\right) p_{1}\left(
s\right) ds\right)  \\
&=&\overline{\theta }_{1}\left( f_{1}\left( u_{n}\left( r\right)
+v_{n}\left( r\right) ,u_{n}\left( r\right) +v_{n}\left( r\right) \right)
\right) P_{1}^{\prime }\left( r\right) ,
\end{eqnarray*}%
and%
\begin{equation*}
v_{n}^{\prime }\left( r\right) \leq \overline{\theta }_{2}\left( f_{2}\left(
u_{n}\left( r\right) +v_{n}\left( r\right) ,u_{n}\left( r\right)
+v_{n}\left( r\right) \right) \right) P_{2}^{\prime }\left( r\right) .
\end{equation*}%
It follows from these last inequalities that 
\begin{equation}
\frac{\left( u_{n}\left( r\right) +v_{n}\left( r\right) \right) ^{\prime }}{%
\left( \overline{\theta }_{1}\left( f_{1}\right) +\overline{\theta }%
_{2}\left( f_{2}\right) \right) \left( \left( u_{n}\left( r\right)
+v_{n}\left( r\right) ,u_{n}\left( r\right) +v_{n}\left( r\right) \right)
\right) }\leq P_{1}^{\prime }\left( r\right) +P_{2}^{\prime }\left( r\right)
,  \label{mat}
\end{equation}%
from which inequality we obtain%
\begin{equation*}
\int_{a_{1}+a_{2}}^{u_{n}\left( r\right) +v_{n}\left( r\right) }\frac{1}{%
\overline{\theta }_{1}\left( f_{1}\left( t,t\right) \right) +\overline{%
\theta }_{2}\left( f_{2}\left( t,t\right) \right) }dt\leq P_{1}\left(
r\right) +P_{2}\left( r\right) .
\end{equation*}%
Now we have 
\begin{equation}
Z\left( u_{n}\left( r\right) +v_{n}\left( r\right) \right) \leq P_{1}\left(
r\right) +P_{2}\left( r\right) ,  \label{zc1}
\end{equation}%
which will play a basic role in the proof of our main results. The
inequalities (\ref{zc1}) can be rewritten as 
\begin{equation}
u_{n}\left( r\right) +v_{n}\left( r\right) \leq Z^{-1}\left( P_{1}\left(
r\right) +P_{2}\left( r\right) \right) .  \label{zc2}
\end{equation}%
This can be easily seen from the fact that $Z$ is a bijection with the
inverse function $Z$ strictly increasing on $\left[ 0,Z\left( \infty \right)
\right) $. Let $M_{1}\geq \max \left\{ 1,\frac{1}{a_{1}}\right\} $ and $%
M_{2}\geq \max \left\{ 1,\frac{1}{a_{2}}\right\} $. The next step is to
integrate (\ref{gen1}) from $0$ to $r$ and bearing in mind (\ref{c22}), we
find 
\begin{eqnarray}
\left( u_{n}\left( r\right) \right) ^{\prime } &\leq &\Psi _{1}^{-1}\left( 
\frac{1}{\xi _{1}\left( r\right) }\int_{0}^{r}\xi _{1}\left( s\right)
p_{1}\left( s\right) f_{1}\left( u_{n}\left( s\right) ,v_{n}\left( s\right)
\right) ds\right)   \notag \\
&\leq &\Psi _{1}^{-1}\left( \frac{1}{\xi _{1}\left( r\right) }%
\int_{0}^{r}\xi _{1}\left( s\right) p_{1}\left( s\right) f_{1}\left(
u_{n}\left( s\right) ,2u_{n}\left( s\right) +v_{n}\left( s\right) \right)
ds\right)   \notag \\
&\leq &\Psi _{1}^{-1}\left( \frac{1}{\xi _{1}\left( r\right) }%
\int_{0}^{r}\xi _{1}\left( s\right) p_{1}\left( s\right) f_{1}\left(
u_{n}\left( s\right) ,u_{n}\left( s\right) +Z^{-1}\left( P_{1}\left(
s\right) +P_{2}\left( s\right) \right) \right) ds\right)   \notag \\
&=&\Psi _{1}^{-1}\left( \frac{1}{\xi _{1}\left( r\right) }\int_{0}^{r}\xi
_{1}\left( s\right) p_{1}\left( s\right) f_{1}\left( u_{n}\left( s\right)
,u_{n}\left( s\right) (1+\frac{1}{u_{n}\left( s\right) }Z^{-1}\left(
P_{1}\left( s\right) +P_{2}\left( s\right) \right) )\right) ds\right)  
\notag \\
&\leq &\Psi _{1}^{-1}\left( \frac{1}{\xi _{1}\left( r\right) }%
\int_{0}^{r}\xi _{1}\left( s\right) p_{1}\left( s\right) f_{1}\left(
u_{n}\left( s\right) ,u_{n}\left( s\right) (1+\frac{1}{a_{1}}Z^{-1}\left(
P_{1}\left( s\right) +P_{2}\left( s\right) \right) )\right) ds\right) 
\label{exin} \\
&\leq &\Psi _{1}^{-1}\left( \frac{1}{\xi _{1}\left( r\right) }%
\int_{0}^{r}\xi _{1}\left( s\right) p_{1}\left( s\right) f_{1}\left(
u_{n}\left( s\right) ,M_{1}\left( 1+Z^{-1}\left( P_{1}\left( s\right)
+P_{2}\left( s\right) \right) \right) \right) ds\right)   \notag \\
&\leq &\Psi _{1}^{-1}\left( \frac{1}{\xi _{1}\left( r\right) }h_{1}\left(
u_{n}\left( r\right) ,M_{1}u_{n}\left( r\right) \right) \int_{0}^{r}\xi
_{1}\left( s\right) p_{1}\left( s\right) \overline{f}_{1}\left(
1+Z^{-1}\left( P_{1}\left( s\right) +P_{2}\left( s\right) \right) \right)
ds\right)   \notag \\
&\leq &\overline{\theta }_{1}\left( h_{1}\left( u_{n}\left( r\right)
,M_{1}u_{n}\left( r\right) \right) \right) \Psi _{1}^{-1}\left( \frac{1}{\xi
_{1}\left( r\right) }\int_{0}^{r}\xi _{1}\left( s\right) p_{1}\left(
s\right) \overline{f}_{1}\left( 1+Z^{-1}\left( P_{1}\left( s\right)
+P_{2}\left( s\right) \right) \right) ds\right)   \notag \\
&=&\overline{\theta }_{1}\left( h_{1}\left( u_{n}\left( r\right)
,M_{1}u_{n}\left( r\right) \right) \right) \overline{P}_{1}^{\prime }\left(
r\right) .  \notag
\end{eqnarray}%
Dividing the inequality (\ref{exin}) by $\overline{\theta }_{1}\left(
h_{1}\left( u_{n}\left( r\right) ,M_{1}u_{n}\left( r\right) \right) \right) $
we see that 
\begin{equation}
\frac{\left( u_{n}\left( r\right) \right) ^{\prime }}{\overline{\theta }%
_{1}\left( h_{1}\left( u_{n}\left( r\right) ,M_{1}u_{n}\left( r\right)
\right) \right) }\leq \overline{P}_{1}^{\prime }\left( r\right) .
\label{mat2}
\end{equation}%
Integrate (\ref{mat2}), from $0$ to $r$, we have%
\begin{equation*}
\int_{a_{1}}^{u_{n}\left( r\right) }\frac{1}{\overline{\theta }_{1}\left(
h_{1}\left( t,M_{1}t\right) \right) }dt\leq \overline{P}_{1}\left( r\right) ,%
\text{ }
\end{equation*}%
which is the same as%
\begin{equation}
\mathcal{H}_{1}\left( u_{n}\left( r\right) \right) \leq \overline{P}%
_{1}\left( r\right) .  \label{ints}
\end{equation}%
Now, we can easy see that $\mathcal{H}_{1}$ is a bijection with the inverse
function $\mathcal{H}_{1}^{-1}$ strictly increasing on $\left[ 0,\mathcal{H}%
_{1}\left( \infty \right) \right) $. By combining this with the previous
inequality, leads to%
\begin{equation}
u_{n}\left( r\right) \leq \mathcal{H}_{1}^{-1}\left( \overline{P}_{1}\left(
r\right) \right) .  \label{int}
\end{equation}%
Returning to $\left\{ v_{n}\left( r\right) \right\} _{n\geq 0}$, one can
show that%
\begin{equation*}
\left( v_{n}\left( r\right) \right) ^{\prime }\leq \overline{\theta }%
_{2}\left( h_{2}\left( v_{n}\left( r\right) ,M_{2}v_{n}\left( r\right)
\right) \right) \overline{P}_{2}^{\prime }\left( r\right) .
\end{equation*}%
Integrating this ordinary differential inequality we get%
\begin{equation*}
\mathcal{H}_{2}\left( v_{n}\left( r\right) \right)
=\int_{a_{2}}^{v_{n}\left( r\right) }\frac{1}{\overline{\theta }_{2}\left(
h_{2}\left( t,M_{2}t\right) \right) }dt\leq \overline{P}_{2}\left( r\right) .
\end{equation*}%
From this inequality, we derive 
\begin{equation}
v_{n}\left( r\right) \leq \mathcal{H}_{2}^{-1}\left( \overline{P}_{2}\left(
r\right) \right) .  \label{int2}
\end{equation}%
In summary, we have found upper bounds for $\left\{ u_{n}\right\} _{n\geq 0}$
and $\left\{ v_{n}\right\} _{n\geq 0}$ which are dependent of $r$. Now let
us complete the proof of Theorem \ref{th1}. We prove that the sequences $%
\left\{ u_{n}\right\} _{n\geq 0}$ and $\left\{ v_{n}\right\} _{n\geq 0}$ are
bounded and equicontinuous on $\left[ 0,c_{0}\right] $ for arbitrary $c_{0}>0
$. Indeed, since 
\begin{equation*}
\left( u_{n}\left( r\right) \right) ^{^{\prime }}\geq 0\text{ and }\left(
v_{n}\left( r\right) \right) ^{^{\prime }}\geq 0\text{ for all }r\geq 0,
\end{equation*}%
it follows that 
\begin{equation*}
u_{n}\left( r\right) \leq u_{n}\left( c_{0}\right) \leq C_{1}\text{ and }%
v_{n}\left( r\right) \leq v_{n}\left( c_{0}\right) \leq C_{2}\text{ on }%
\left[ 0,c_{0}\right] .
\end{equation*}%
Here $C_{1}=\mathcal{H}_{1}^{-1}\left( \overline{P}_{1}\left( c_{0}\right)
\right) $ and $C_{2}=\mathcal{H}_{2}^{-1}\left( \overline{P}_{2}\left(
c_{0}\right) \right) $ are positive constants. Recall that $\left\{
u_{n}\right\} _{n\geq 0}$ and $\left\{ v_{n}\right\} _{n\geq 0}$ are bounded
on $\left[ 0,c_{0}\right] $ for arbitrary $c_{0}>0$. Using this fact, we
show that the same is true of $\left( u_{n}\left( r\right) \right) ^{\prime }
$ and $\left( v_{n}\left( r\right) \right) ^{\prime }$. Indeed, for any $%
r\geq 0$,%
\begin{eqnarray*}
\left( u_{n}\left( r\right) \right) ^{\prime } &=&\Psi _{1}^{-1}\left( \frac{%
1}{\xi _{1}\left( r\right) }\int_{0}^{r}\xi _{1}\left( s\right) p_{1}\left(
s\right) f_{1}\left( u_{n-1}\left( s\right) ,v_{n-1}\left( s\right) \right)
ds\right)  \\
&\leq &\Psi _{1}^{-1}\left( \frac{1}{\xi _{1}\left( r\right) }%
\int_{0}^{r}\xi _{1}\left( s\right) p_{1}\left( s\right) f_{1}\left(
u_{n}\left( s\right) ,v_{n}\left( s\right) \right) ds\right)  \\
&\leq &\Psi _{1}^{-1}\left( \left\Vert p_{1}\right\Vert _{\infty
}f_{1}\left( C_{1},C_{2}\right) \frac{1}{\xi _{1}\left( r\right) }%
\int_{0}^{r}\xi _{1}\left( s\right) ds\right)  \\
&\leq &\Psi _{1}^{-1}\left( \left\Vert p_{1}\right\Vert _{\infty
}f_{1}\left( C_{1},C_{2}\right) \int_{0}^{r}ds\right)  \\
&\leq &\Psi _{1}^{-1}\left( \left\Vert p_{1}\right\Vert _{\infty
}f_{1}\left( C_{1},C_{2}\right) c_{0}\right) \text{ on }\left[ 0,c_{0}\right]
.
\end{eqnarray*}%
Similar arguments show that%
\begin{eqnarray*}
\left( v_{n}\left( r\right) \right) ^{\prime } &=&\Psi _{2}^{-1}\left( \frac{%
1}{\xi _{2}\left( r\right) }\int_{0}^{r}\xi _{2}\left( s\right) p_{2}\left(
s\right) f_{2}\left( u_{n-1}\left( s\right) ,v_{n-1}\left( s\right) \right)
ds\right)  \\
&\leq &\Psi _{2}^{-1}\left( \left\Vert p_{2}\right\Vert _{\infty
}f_{2}\left( C_{1},C_{2}\right) c_{0}\right) \text{ on }\left[ 0,c_{0}\right]
.
\end{eqnarray*}%
It remains, to prove that $\left\{ u_{n}\right\} _{n\geq 0}$ and $\left\{
v_{n}\right\} _{n\geq 0}$ are equicontinuous on $\left[ 0,c_{0}\right] $ for
arbitrary $c_{0}>0$. Let $\varepsilon _{1}$, $\varepsilon _{2}>0$. To verify
equicontinuous on $\left[ 0,c_{0}\right] $, observe that%
\begin{eqnarray*}
\left\vert u_{n}\left( x\right) -u_{n}\left( y\right) \right\vert 
&=&\left\vert \left( u_{n}\left( \xi _{1}\right) \right) ^{\prime
}\right\vert \left\vert x-y\right\vert \leq \Psi _{1}^{-1}\left( \left\Vert
p_{1}\right\Vert _{\infty }f_{1}\left( C_{1},C_{2}\right) c_{0}\right)
\left\vert x-y\right\vert , \\
\left\vert v_{n}\left( x\right) -v_{n}\left( y\right) \right\vert 
&=&\left\vert \left( v_{n}\left( \xi _{2}\right) \right) ^{\prime
}\right\vert \left\vert x-y\right\vert \leq \Psi _{2}^{-1}\left( \left\Vert
p_{2}\right\Vert _{\infty }f_{2}\left( C_{1},C_{2}\right) c_{0}\right)
\left\vert x-y\right\vert ,
\end{eqnarray*}%
for all $n\in \mathbb{N}$ and all $x,y\in \left[ 0,c_{0}\right] $ and for $%
\xi _{1}$, $\xi _{2}$ the constants from the mean value theorem. So it
suffices to take%
\begin{equation*}
\delta _{1}=\frac{\varepsilon _{1}}{\Psi _{1}^{-1}\left( \left\Vert
p_{1}\right\Vert _{\infty }f_{1}\left( C_{1},C_{2}\right) c_{0}\right) }%
\text{ and }\delta _{2}=\frac{\varepsilon _{2}}{\Psi _{2}^{-1}\left(
\left\Vert p_{2}\right\Vert _{\infty }f_{2}\left( C_{1},C_{2}\right)
c_{0}\right) },
\end{equation*}%
to see that $\left\{ u_{n}\right\} _{n\geq 0}$ and $\left\{ v_{n}\right\}
_{n\geq 0}$ are equicontinuous on $\left[ 0,c_{0}\right] $. In particular,
it follows from the Arzela--Ascoli theorem that there exists a function $%
u\in C\left( \left[ 0,c_{0}\right] \right) $ and a subsequence $N_{1}$ of $%
\mathbb{N}^{\ast }$ with $u_{n}\left( r\right) $ converging uniformly on $u$
to $\left[ 0,c_{0}\right] $ as $n\rightarrow \infty $ through $N_{1}$. By
the same token there exists a function $v\in C\left( \left[ 0,c_{0}\right]
\right) $ and a subsequence $N_{2}$ of $\mathbb{N}^{\ast }$ with $%
v_{n}\left( r\right) $ converging uniformly to $v$ on $\left[ 0,c_{0}\right] 
$ as $n\rightarrow \infty $ through $N_{2}$. Thus $\left\{ \left(
u_{n}\left( r\right) ,v_{n}\left( r\right) \right) \right\} _{n\in N_{2}}$
converges uniformly on $\left[ 0,c_{0}\right] $ to $\left( u,v\right) \in
C\left( \left[ 0,c_{0}\right] \right) \times C\left( \left[ 0,c_{0}\right]
\right) $ through $N_{2}$ (see L\"{u}-O'Regan-Agarwal \cite{LU}). The limit
function $\left( u,v\right) $ constructed in this way will be nonnegative,
radially symmetric and nondecreasing with respect to $r$ and is a solution
of system (\ref{11}). Moreover, the radial solutions of (\ref{11}) with $%
u\left( 0\right) =a_{1},$ $v\left( 0\right) =a_{2}$ satisfy:%
\begin{eqnarray}
u\left( r\right)  &=&a_{1}+\int_{0}^{r}\Psi _{1}^{-1}\left( \frac{1}{\xi
_{1}\left( t\right) }\int_{0}^{t}\xi _{1}\left( s\right) p_{1}\left(
s\right) f_{1}\left( u\left( s\right) ,v\left( s\right) \right) ds\right) dt,%
\text{ }r\geq 0,  \label{eq1} \\
v\left( r\right)  &=&a_{2}+\int_{0}^{r}\Psi _{2}^{-1}\left( \frac{1}{\xi
_{2}\left( t\right) }\int_{0}^{t}\xi _{2}\left( s\right) p_{2}\left(
s\right) f_{2}\left( u\left( s\right) ,v\left( s\right) \right) ds\right) dt,%
\text{ }r\geq 0.  \label{eq2}
\end{eqnarray}%
In the case $\underline{P}_{1}\left( \infty \right) =\underline{P}_{2}\left(
\infty \right) =\infty $, we observe that 
\begin{eqnarray}
u\left( r\right)  &=&a_{1}+\int_{0}^{r}\Psi _{1}^{-1}\left( \frac{1}{\xi
_{1}\left( t\right) }\int_{0}^{t}\xi _{1}\left( s\right) p_{1}\left(
s\right) f_{1}\left( u\left( s\right) ,v\left( s\right) \right) ds\right) dt
\notag \\
&\geq &a_{1}+\int_{0}^{r}\Psi _{1}^{-1}\left( \frac{1}{\xi _{1}\left(
t\right) }\int_{0}^{t}\xi _{1}\left( s\right) p_{1}\left( s\right)
f_{1}\left( a_{1},a_{2}+\underline{\theta }_{2}(f_{2}\left(
a_{1},a_{2}\right) )P_{2}\left( s\right) \right) ds\right) dt  \label{ints1}
\\
&=&a_{1}+\underline{P}_{1}\left( r\right) .  \notag
\end{eqnarray}%
We repeat the argument applied in the proof of (\ref{ints1}) 
\begin{eqnarray}
v\left( r\right)  &=&a_{2}+\int_{0}^{r}\Psi _{2}^{-1}\left( \frac{1}{\xi
_{2}\left( t\right) }\int_{0}^{t}\xi _{2}\left( s\right) p_{2}\left(
s\right) f_{2}\left( u\left( s\right) ,v\left( s\right) \right) ds\right) dt
\notag \\
&\geq &a_{2}+\int_{0}^{r}\Psi _{2}^{-1}\left( \frac{1}{\xi _{2}\left(
t\right) }\int_{0}^{t}\xi _{2}\left( s\right) p_{2}\left( s\right)
f_{2}\left( a_{1}+\underline{\theta }_{1}(f_{2}\left( a_{1},a_{2}\right)
)P_{1}\left( s\right) ,a_{2}\right) ds\right) dt  \label{ints2} \\
&=&a_{2}+\underline{P}_{2}\left( r\right) .  \notag
\end{eqnarray}%
By taking limits in (\ref{ints1}) and (\ref{ints2}), we get entire large
solutions 
\begin{equation*}
\lim_{r\rightarrow \infty }u\left( r\right) =\infty \text{ and }%
\lim_{r\rightarrow \infty }v\left( r\right) =\infty .
\end{equation*}%
Consequently, $\left( u,v\right) $ is an entire large solution of (\ref{11}).

The next purpose of the paper is to give a sufficient condition to obtain an
entire bounded solution to (\ref{11}). Our result in this case is the
following:

\begin{theorem}
\label{th2}The system (\ref{11}) has one positive radial solution $\left(
u,v\right) \in C^{1}\left( \left[ 0,\infty \right) \right) \times
C^{1}\left( \left[ 0,\infty \right) \right) $ given that\textit{\ }$\mathcal{%
H}_{1}\left( \infty \right) =\mathcal{H}_{2}\left( \infty \right) =\infty $
and \textrm{(P1)}, \textrm{(C1)}, \textrm{(C2)} hold true. Moreover, if $%
\overline{P}_{1}\left( \infty \right) <\infty $ and $\overline{P}_{2}\left(
\infty \right) <\infty $ then 
\begin{equation*}
\lim_{r\rightarrow \infty }u\left( r\right) <\infty \text{ and }%
\lim_{r\rightarrow \infty }v\left( r\right) <\infty .
\end{equation*}
\end{theorem}

\subparagraph{\textbf{Proof of Theorem \protect\ref{th2}:}}

The existence part is proved in Theorem \ref{th1}. Assume $\overline{P}%
_{1}\left( \infty \right) <\infty $ and $\overline{P}_{2}\left( \infty
\right) <\infty $. Proceeding as in the proof of (\ref{int}) and (\ref{int2}%
) with the integral equations (\ref{eq1}) and (\ref{eq2}), one gets the
estimates%
\begin{equation*}
u\left( r\right) \leq \mathcal{H}_{1}^{-1}\left( \overline{P}_{1}\left(
\infty \right) \right) <\infty \text{ and }v\left( r\right) \leq \mathcal{H}%
_{2}^{-1}\left( \overline{P}_{2}\left( \infty \right) \right) <\infty \text{
for all }r\geq 0.
\end{equation*}%
Thus $\left( u,v\right) $ is a positive entire bounded solution of the
system (\ref{11}).

Concerning the existence of semifinite entire large solutions to (\ref{11}),
we have the following:

\begin{theorem}
\label{th3}The system (\ref{11}) has one positive radial solution $\left(
u,v\right) \in C^{1}\left( \left[ 0,\infty \right) \right) \times
C^{1}\left( \left[ 0,\infty \right) \right) $ given that $\mathcal{H}%
_{1}\left( \infty \right) =\mathcal{H}_{2}\left( \infty \right) =\infty $
and \textrm{(P1)}, \textrm{(C1)}, \textrm{(C2)} hold true. Moreover, the
following hold:

1)\quad If $\overline{P}_{1}\left( \infty \right) <\infty $ and $\underline{P%
}_{2}\left( \infty \right) =\infty $ then 
\begin{equation*}
\lim_{r\rightarrow \infty }u\left( r\right) <\infty \text{ and }%
\lim_{r\rightarrow \infty }v\left( r\right) =\infty .
\end{equation*}

2)\quad If $\underline{P}_{1}\left( \infty \right) =\infty $ and $\overline{P%
}_{2}\left( \infty \right) <\infty $ then 
\begin{equation*}
\lim_{r\rightarrow \infty }u\left( r\right) =\infty \text{ and }%
\lim_{r\rightarrow \infty }v\left( r\right) <\infty .
\end{equation*}
\end{theorem}

\subparagraph{\textbf{Proof of Theorem \protect\ref{th3}:}}

The existence part is proved in Theorem \ref{th1}.

\textbf{1):} As in the proof of Theorem \ref{th1} and Theorem \ref{th2}, we
have 
\begin{equation*}
u\left( r\right) \leq \mathcal{H}_{1}^{-1}\left( \overline{P}_{1}\left(
\infty \right) \right) <\infty \text{ and }v\left( r\right) \geq a_{2}+%
\underline{P}_{1}\left( r\right) .
\end{equation*}%
Observing that $\overline{P}_{1}\left( \infty \right) <\infty $ and $%
\underline{P}_{2}\left( \infty \right) =\infty $ the above relations yield 
\begin{equation*}
\lim_{r\rightarrow \infty }u\left( r\right) <\infty \text{ and }%
\lim_{r\rightarrow \infty }v\left( r\right) =\infty \text{.}
\end{equation*}%
This completes the proof.

\textbf{2): }Arguing as above, we obtain%
\begin{equation}
u\left( r\right) \geq a_{1}+\underline{P}_{1}\left( r\right) \text{ and }%
v\left( r\right) \leq \mathcal{H}_{2}^{-1}\left( \overline{P}_{2}\left(
r\right) \right) .  \label{t2}
\end{equation}%
Our conclusion follows now by letting $r\rightarrow \infty $ in (\ref{t2}).

We now propose a more refined question concerning the solutions of system (%
\ref{11}).\ In analogy with Theorems \ref{th1}-\ref{th3}, we can also prove
the following three theorems. The first is the following:

\begin{theorem}
\label{th4}The system (\ref{11}) has one positive bounded radial solution $%
\left( u,v\right) \in C^{1}\left( \left[ 0,\infty \right) \right) \times
C^{1}\left( \left[ 0,\infty \right) \right) $ given that $\overline{P}%
_{1}\left( \infty \right) <\mathcal{H}_{1}\left( \infty \right) <\infty $, $%
\overline{P}_{2}\left( \infty \right) <\mathcal{H}_{2}\left( \infty \right)
<\infty $, \textrm{(P1), (C1),\ (C2) }hold true. \ Moreover,%
\begin{equation*}
\left\{ 
\begin{array}{l}
a_{1}+\underline{P}_{1}\left( r\right) \leq u\left( r\right) \leq \mathcal{H}%
_{1}^{-1}\left( \overline{P}_{1}\left( r\right) \right) , \\ 
a_{2}+\underline{P}_{1}\left( r\right) \leq v\left( r\right) \leq \mathcal{H}%
_{2}^{-1}\left( \overline{P}_{2}\left( r\right) \right) .%
\end{array}%
\right. \text{ }
\end{equation*}
\end{theorem}

\subparagraph{\textbf{Proof of Theorem \protect\ref{th4}: }}

The existence part is proved in Theorem \ref{th1}. Next, by a simple
calculation together with (\ref{ints}) and the conditions of the theorem we
obtain:%
\begin{equation*}
\mathcal{H}_{1}\left( u_{n}\left( r\right) \right) \leq \overline{P}%
_{1}\left( \infty \right) <\mathcal{H}_{1}\left( \infty \right) <\infty 
\text{ and }v_{n}\left( r\right) \leq \mathcal{H}_{2}^{-1}\left( \overline{P}%
_{2}\left( \infty \right) \right) <\infty .
\end{equation*}%
On the other hand, since $\mathcal{H}_{1}^{-1}$ is strictly increasing on $%
\left[ 0,\mathcal{H}_{1}\left( \infty \right) \right) $, we find that 
\begin{equation*}
u_{n}\left( r\right) \leq \mathcal{H}_{1}^{-1}\left( \overline{P}_{1}\left(
\infty \right) \right) <\infty ,
\end{equation*}%
and then the non-decreasing sequences $\left\{ u_{n}\right\} _{n\geq 0}$ and 
$\left\{ v_{n}\right\} _{n\geq 0}$ are bounded above for all $r\geq 0$ and
all $n$. Now we use this observation to conclude 
\begin{equation*}
\left( u_{n}\left( r\right) ,v_{n}\left( r\right) \right) \overset{%
n\rightarrow \infty }{\rightarrow }\left( u\left( r\right) ,v\left( r\right)
\right)
\end{equation*}%
and then the limit functions $u$ and $v$ are positive entire bounded radial
solutions of system (\ref{11}).\textbf{\ }This completes the proof.

\begin{theorem}
\label{th5}Assume \textrm{(P1), (C1) }and\textrm{\ (C2) }hold true. The
following hold true:

i)\quad The system (\ref{11}) has one positive radial solution $\left(
u,v\right) \in C^{1}\left( \left[ 0,\infty \right) \right) \times
C^{1}\left( \left[ 0,\infty \right) \right) $ such that $\lim_{r\rightarrow
\infty }u\left( r\right) =\infty $ and $\lim_{r\rightarrow \infty }v\left(
r\right) <\infty $ given that $\mathcal{H}_{1}\left( \infty \right) =\infty $%
, $\underline{P}_{1}\left( \infty \right) =\infty $ and $\overline{P}%
_{2}\left( \infty \right) <\mathcal{H}_{2}\left( \infty \right) <\infty $.

ii)\quad The system (\ref{11}) has one positive radial solution $\left(
u,v\right) \in C^{1}\left( \left[ 0,\infty \right) \right) \times
C^{1}\left( \left[ 0,\infty \right) \right) $ such that $\lim_{r\rightarrow
\infty }u\left( r\right) <\infty $ and $\lim_{r\rightarrow \infty }v\left(
r\right) =\infty $ given that $\overline{P}_{1}\left( \infty \right) <%
\mathcal{H}_{1}\left( \infty \right) <\infty $ and $\mathcal{H}_{2}\left(
\infty \right) =\infty $, $\underline{P}_{2}\left( \infty \right) =\infty $.
\end{theorem}

\subparagraph{\textbf{Proof of Theorem \protect\ref{th5}:}}

The proof for these cases is similar as the above and is therefore omitted.

\end{document}